\documentclass[11pt]{article}
\usepackage[a4paper]{anysize}\marginsize{3.5cm}{3.5cm}{1.3cm}{2cm}
\pdfpagewidth=\paperwidth \pdfpageheight=\paperheight
\usepackage{amsfonts,amssymb,amsthm,amsmath,eucal}
\usepackage{pgf}
\usepackage{bbm}

\usepackage{tikz} 
\usepackage{float}
\usetikzlibrary{arrows}
\usepackage{subfigure}
\usepackage{caption}

\theoremstyle{plain}
\newtheorem{thm}{Theorem}[section]
\newtheorem{theorem}[thm]{Theorem}

\theoremstyle{definition}

\newtheorem{problem}[thm]{Problem}

\newtheorem{thevarthm}[thm]{\varthmname}

\newenvironment{varthm*}[1]{\trivlist\item[]{\bf #1.}\it}{\endtrivlist}

\renewcommand\geq{\geqslant}

\renewcommand\leq{\leqslant}

\newcommand\be{\begin{eqnarray*}}
\newcommand\ee{\end{eqnarray*}}

\newcommand\K{\mathbb K}

\renewcommand\P{\mathbb P}

\newcommand\newop[2]{\def#1{\mathop{\rm #2}\nolimits}}
\newop\log{log}
\newop\ord{ord}
\newop\Gal{Gal}
\newop\SL{SL}
\newop\Bl{Bl}
\newop\mult{mult}
\newop\mass{mass}
\newop\div{div}
\newop\codim{codim}
\newop\sing{sing}
\newop\Zeroes{Zeroes}
\newop\Ass{Ass}
\newop\reg{reg}
\newop\Ext{Ext}
\newop\PSL{PSL}
\newop{\symassreg}{areg}

\newcommand\lra\longrightarrow

\def\keywordname{{\bfseries Keywords}}%
\def\keywords#1{\par\addvspace\medskipamount{\rightskip=0pt plus1cm
\def\and{\ifhmode\unskip\nobreak\fi\ $\cdot$
}\noindent\keywordname\enspace\ignorespaces#1\par}}
\def\subclassname{{\bfseries Mathematics Subject Classification
(2000)}\enspace}
\def\subclass#1{\par\addvspace\medskipamount{\rightskip=0pt plus1cm
\def\and{\ifhmode\unskip\nobreak\fi\ $\cdot$
}\noindent\subclassname\ignorespaces#1\par}}

\definecolor{uuuuuu}{rgb}{0,0,1}
\definecolor{qqqqff}{rgb}{0,0,1}
\definecolor{xdxdff}{rgb}{0,0,1}

\begin{document}

\author{G.~Malara, J.~Szpond}
\title{The containment problem and a rational simplicial arrangement}
\date{\today}
\maketitle
\thispagestyle{empty}

\begin{abstract}
   Since Dumnicki, Szemberg and Tutaj-Gasi\'nska gave in 2013 in \cite{DST13b}
   the first example of a set of points in the complex projective plane such that
   for its homogeneous ideal $I$ the containment of the third symbolic power
   in the second ordinary power fails, there has been considerable interest
   in searching for further examples with this property and investigating into
   the nature of such examples. Many examples, defined over various fields,
   have been found but so far there has been essentially just one example found
   of $19$ points
   defined over the rationals, see \cite[Theorem A, Problem 1]{LBS17}.
   In \cite[Problem 5.1]{FKLT} the authors asked if there are other
   rational examples. This has motivated our research.
   The purpose of this note is to flag the existence of a new example
   of a set of $49$ rational points with the same non-containment property
   for powers of its homogeneous ideal. Here we establish the existence and
   justify it computationally. A more conceptual proof, based on Seceleanu's
   criterion \cite{Sec15} will be published elsewhere \cite{MalSzp}.
\subclass{MSC 13F20 \and MSC 14N20 \and MSC 13A15 \and MSC 52C35}
\end{abstract}


\section{Introduction}
   The begin of the Millennium has brought ground breaking results on the containment
   relation between symbolic and ordinary powers of homogeneous ideals established
   by Ein-Lazarsfeld-Smith \cite{ELS01} in characteristic $0$ and Hochster-Huneke \cite{HoHu02}
   in positive characteristic.

   Let $I$ be a homogeneous ideal in the ring of polynomials $R=\K[x_0,\ldots,x_N]$ defined over a field $\K$.
   For a positive integer $m$ one defines the $m$-th \emph{symbolic power} $I^{(m)}$ of $I$ as
   \begin{equation}\label{eq:symbolic power}
      I^{(m)}=\bigcap_{P\in\Ass(I)}\left(I^mR_P\cap R\right),
   \end{equation}
   where $\Ass(I)$ is the set of associated primes of $I$.
   A lot of research has been motivated in recent two decades by the following central question.
\begin{problem}[The containment problem]\label{pro:containment}
   Decide for which $m$ and $r$ there is the containment
   \begin{equation}\label{eq:containment m in r}
      I^{(m)}\subset I^r.
   \end{equation}
\end{problem}
   Whereas to decide the reverse containment i.e. to decide when
   $$I^r\subset I^{(m)}$$
   holds is almost elementary: It does
   if and only if $m\leq r$ \cite[Lemma 8.4.1]{PSC},
   the elegant uniform answer to Problem \ref{pro:containment},
   provided in the following Theorem, came as a surprise.
\begin{theorem}[Ein-Lazarsfeld-Smith, Hochster-Huneke]\label{thm:ELS}
   Let $I\subset\K[x_0,\ldots,x_N]$ be a homogeneous ideal such that every
   component of its zero locus $V(I)$ has codimension at most $e$. Then the containment
   $$I^{(m)}\subset I^r$$
   holds for all $m\geq er$.
\end{theorem}
   It is natural to wonder to what extend the bound in Theorem \ref{thm:ELS}
   is optimal. A lot of attention has been given to this problem in general.
   Here we focus on the simples non-trivial situation. If $I$ is an ideal
   of points in $\P^2$, then Theorem \ref{thm:ELS} asserts that there is always
   the containment
   $$I^{(4)}\subset I^2.$$
   It is also elementary to provide examples when the containment $I^{(2)}\subset I^2$
   fails (three general points in $\P^2$ are sufficient). The question if in the remaining case,
   i.e.
   \begin{equation}\label{eq:I3inI2}
      I^{(3)}\subset I^2
   \end{equation}
   the containment always holds has been raised by Huneke and repeated in the literature
   by several authors who proved a considerable number of special cases. However,
   in 2013 Dumnicki, Szemberg and Tutaj-Gasi\'nska showed in \cite{DST13b} that
   the containment in \eqref{eq:I3inI2} fails for an explicit set of $12$ points in the complex
   projective plane, the points dual to the Hesse arrangement of lines, see \cite{ArtDol09}.
   These points arise also as intersection points of certain arrangement of $9$ lines.
   These lines intersect by $3$ in configuration points and there are no other intersection
   points among them. This the only known non-trivial arrangement of lines in characteristic $0$
   where each intersection point of a pair of lines lies on exactly three distinct lines.

   Soon after \cite{DST13b} became public, series of additional examples in all characteristics
   have been constructed, see e.g. \cite{BNAL}, \cite{Real}, \cite{HarSec15}, \cite{LBS17}, \cite{NagSec16}
   and \cite{SzeSzp17} for an overview.
   All these constructions were based on the same idea of checking the containment for the set
   of intersection points of multiplicity at least $3$ (at least $3$ lines pass through a point)
   of some special arrangements of lines. The relevance of points of multiplicity at least $3$
   follows from the Zariski-Nagata Theorem, see \cite[Theorem 3.14]{Eisenbud}, which, at least
   in characteristic $0$ characterizes the $m$-th symbolic power of a radical ideal $I$
   as the set of polynomials vanishing to order at least $m$ along the set of zeroes of $I$.
   Among these arrangements prominent role is played by arrangements coming from finite
   reflection groups, see e.g. \cite{NC}. Finding configurations of rational points turned
   out to be a little bit harder. The only known example of $19$ points arising as intersection
   points of multiplicity at least $3$ of an arrangement of $12$ lines is described in
   the recent work of Lampa-Baczy\'nska and the second author. In this note we exhibit an
   additional example.

\section{The simplicial arrangement $A(25,2)$}
   An arrangement of lines is a complex generated in the real projective plane by a finite
   family of lines that do not form a pencil. An arrangement is \emph{simplicial} if all its
   faces (connected components of the complement of the union of lines) are triangles.
   Apparently, simplicial arrangements have been first studied by Melchior in connection
   with his proof of the celebrated Sylvester-Gallai Theorem, see \cite{Mel41}.
   They have received considerable attention after Deligne's paper \cite{Del72}.
   Recently Gr\"unbaum \cite{Gru09} classified all known simplicial arrangement.
   Cuntz pointed out that there might be gaps in Gr\"unbaum's list, constructed
   $4$ new examples and provided a conceptual proof of the completeness of the classification
   for arrangements of up to $27$ lines. The classification in general seems to be an open
   and challenging problem.

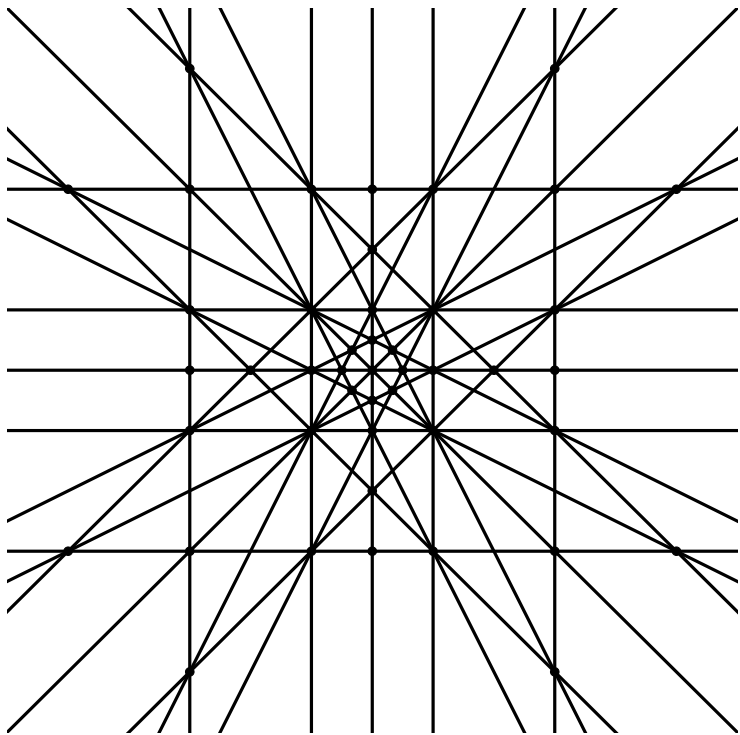
\begin{figure}[H]
\centering
\definecolor{uuuuuu}{rgb}{1,1,1}
\begin{tikzpicture}[line cap=round,line join=round,>=triangle 45,x=1.0cm,y=1.0cm,scale=0.8]
\clip(-6,-6) rectangle (6,6);
\draw [line width=1.2pt,domain=-17.344264336523096:18.050613460530244] plot(\x,{(--6.-0.*\x)/6.});
\draw [line width=1.2pt,domain=-17.344264336523096:18.050613460530244] plot(\x,{(-6.-0.*\x)/6.});
\draw [line width=1.2pt] (-1.,-8.384373662511335) -- (-1.,9.606934147614442);
\draw [line width=1.2pt] (1.,-8.384373662511335) -- (1.,9.606934147614442);
\draw [line width=1.2pt,domain=-17.344264336523096:18.050613460530244] plot(\x,{(--18.-0.*\x)/6.});
\draw [line width=1.2pt,domain=-17.344264336523096:18.050613460530244] plot(\x,{(-18.-0.*\x)/6.});
\draw [line width=1.2pt] (3.,-8.384373662511335) -- (3.,9.606934147614442);
\draw [line width=1.2pt] (-3.,-8.384373662511335) -- (-3.,9.606934147614442);
\draw [line width=1.2pt,domain=-17.344264336523096:18.050613460530244] plot(\x,{(-0.-0.*\x)/6.});
\draw [line width=1.2pt] (0.,-8.384373662511335) -- (0.,9.606934147614442);
\draw [line width=1.2pt,domain=-17.344264336523096:18.050613460530244] plot(\x,{(--8.--4.*\x)/4.});
\draw [line width=1.2pt,domain=-17.344264336523096:18.050613460530244] plot(\x,{(-8.--4.*\x)/4.});
\draw [line width=1.2pt,domain=-17.344264336523096:18.050613460530244] plot(\x,{(-8.--4.*\x)/-4.});
\draw [line width=1.2pt,domain=-17.344264336523096:18.050613460530244] plot(\x,{(-2.-4.*\x)/-2.});
\draw [line width=1.2pt,domain=-17.344264336523096:18.050613460530244] plot(\x,{(--2.-4.*\x)/2.});
\draw [line width=1.2pt,domain=-17.344264336523096:18.050613460530244] plot(\x,{(-2.--2.*\x)/-4.});
\draw [line width=1.2pt,domain=-17.344264336523096:18.050613460530244] plot(\x,{(--2.-2.*\x)/-4.});
\draw [line width=1.2pt,domain=-17.344264336523096:18.050613460530244] plot(\x,{(--2.--4.*\x)/-2.});
\draw [line width=1.2pt,domain=-17.344264336523096:18.050613460530244] plot(\x,{(-2.--4.*\x)/2.});
\draw [line width=1.2pt,domain=-17.344264336523096:18.050613460530244] plot(\x,{(--2.--2.*\x)/4.});
\draw [line width=1.2pt,domain=-17.344264336523096:18.050613460530244] plot(\x,{(-2.-2.*\x)/4.});
\draw [line width=1.2pt,domain=-17.344264336523096:18.050613460530244] plot(\x,{(-0.-6.*\x)/6.});
\draw [line width=1.2pt,domain=-17.344264336523096:18.050613460530244] plot(\x,{(-0.--6.*\x)/6.});
\draw [line width=1.2pt,domain=-17.344264336523096:18.050613460530244] plot(\x,{(-8.-4.*\x)/4.});
\begin{scriptsize}
\draw [fill=black] (0.,0.) circle (2.0pt);
\draw [fill=black] (-3.,-3.) circle (2.0pt);
\draw [fill=black] (3.,-3.) circle (2.0pt);
\draw [fill=black] (3.,3.) circle (2.0pt);
\draw [fill=black] (-3.,3.) circle (2.0pt);
\draw [fill=black] (-3.,1.) circle (2.0pt);
\draw [fill=black] (-3.,-1.) circle (2.0pt);
\draw [fill=black] (-1.,-3.) circle (2.0pt);
\draw [fill=black] (1.,-3.) circle (2.0pt);
\draw [fill=black] (3.,-1.) circle (2.0pt);
\draw [fill=black] (3.,1.) circle (2.0pt);
\draw [fill=black] (1.,3.) circle (2.0pt);
\draw [fill=black] (-1.,3.) circle (2.0pt);
\draw [fill=black] (-1.,1.) circle (2.0pt);
\draw [fill=black] (1.,1.) circle (2.0pt);
\draw [fill=black] (1.,-1.) circle (2.0pt);
\draw [fill=black] (-1.,-1.) circle (2.0pt);
\draw [fill=black] (-3.,0.) circle (2.0pt);
\draw [fill=black] (0.,-3.) circle (2.0pt);
\draw [fill=black] (3.,0.) circle (2.0pt);
\draw [fill=black] (0.,3.) circle (2.0pt);
\draw [fill=black] (-2.,0.) circle (2.0pt);
\draw [fill=black] (2.,0.) circle (2.0pt);
\draw [fill=black] (0.,-2.) circle (2.0pt);
\draw [fill=black] (0.,2.) circle (2.0pt);
\draw [fill=black] (0.,1.) circle (2.0pt);
\draw [fill=black] (-1.,0.) circle (2.0pt);
\draw [fill=black] (0.,-1.) circle (2.0pt);
\draw [fill=black] (1.,0.) circle (2.0pt);
\draw [fill=black] (0.3333333333333333,-0.3333333333333333) circle (2.0pt);
\draw [fill=black] (0.,-0.5) circle (2.0pt);
\draw [fill=black] (-0.3333333333333333,-0.3333333333333333) circle (2.0pt);
\draw [fill=black] (-0.5,0.) circle (2.0pt);
\draw [fill=black] (-0.3333333333333333,0.3333333333333333) circle (2.0pt);
\draw [fill=black] (0.,0.5) circle (2.0pt);
\draw [fill=black] (0.3333333333333333,0.3333333333333333) circle (2.0pt);
\draw [fill=black] (0.5,0.) circle (2.0pt);
\draw [fill=black] (-3.,-5.) circle (2.0pt);
\draw [fill=black] (-5.,-3.) circle (2.0pt);
\draw [fill=black] (3.,-5.) circle (2.0pt);
\draw [fill=black] (5.,-3.) circle (2.0pt);
\draw [fill=black] (5.,3.) circle (2.0pt);
\draw [fill=black] (3.,5.) circle (2.0pt);
\draw [fill=black] (-3.,5.) circle (2.0pt);
\draw [fill=black] (-5.,3.) circle (2.0pt);
\end{scriptsize}
\end{tikzpicture}
\caption{Affine part of the simplicial arrangement $A(25,2)$}
\label{fig:A21}
\end{figure}

   In the present paper we study the arrangement denoted as $A(25,2)$ in Gr\"unbaum's list
   \cite[p. 20]{Gru09}. The same notation is used by Cuntz \cite[p. 699]{Cun12}. This arrangement
   is defined by $25$ lines. They intersect altogether in $85$ points of which $49$ have
   multiplicity at least $3$. This is the set $Z$ of these $49$ points that we are interested
   in. We denote the ideal of $Z$ by $I$. The arrangement can be defined in the rational projective plane and exhibits symmetries
   of a square. In Figure \ref{fig:A21} there are $24$ lines from the $A(25,2)$ arrangement
   visualized. An invisible line is the line at infinity. The product of the equations
   of all $25$ lines defines
   an element in $I^{(3)}$ which is not contained in $I^2$, hence the containment in \eqref{eq:I3inI2} fails for $Z$.
   There are $41$ points from $49$ points in $Z$ visible in Figure \ref{fig:A21}. The
   remaining $8$ points belong to the line at infinity.
   Below we list coordinates of all points in $Z$. In order to alleviate notation
   we use the $\pm$ notation, which indicates
   that points all possible signs should be considered.
   \begin{equation*}
   \begin{array}{ccccc}
      (0:0:1), & (0:\pm 12:1), & (\pm 12: 0:1), & (\pm 6:\pm 6:1), & (\pm 6: \pm 18:1),  \\
      (\pm 3:0:1), & (0:\pm 3:1), & (\pm 6:0:1), & (0:\pm 6:1), & (\pm 2:\pm 2:1), \\
      (\pm 18:\pm 6:1), & (\pm 18: \pm 18:1), & (\pm 18: \pm 30:1) & (\pm 30: \pm 18: 1). &
   \end{array}
   \end{equation*}
   Points at infinity are:
   \begin{equation*}
   \begin{array}{cccc}
      (1:0:0),  & (0:1:0), & (1:-1:0)), & (1:1:0), \\
      (1:-2:0), & (1:2:0), & (-2:1:0), & (2:1:0).
   \end{array}
   \end{equation*}
   In order to list the arrangement lines we use the same $\pm$ convention. The line $z=0$ is the line at infinity.
   \begin{equation*}
   \begin{array}{cccccccc}
      x, & y, & z, & x\pm y, & x\pm 6z, & y\pm 6z, & x\pm 18z, & y\pm 18z,
   \end{array}
   \end{equation*}
   \begin{equation*}
   \begin{array}{ccc}
      y\pm 2x\pm 6z, & y\pm x\pm 12z, & 2y\pm x\pm 6z.
   \end{array}
   \end{equation*}
   It is easy to check with a symbolic algebra program, we used Singular \cite{Singular},
   that $I$ is generated by $3$ reducible octics:
   \begin{equation*}
   \begin{array}{l}
      f:=yz (3553x^6-15102x^4y^2+11549x^2y^4+y^6+385x^4z^2-560x^2y^2z^2\\
      \rule{1cm}{0cm} -189y^4z^2+6300x^2z^4+6804y^2z^4-14665z^6),\\
      g:=xz (-11426x^6+4001x^4y^2-4002x^2y^4+y^6-15874x^4z^2+360x^2y^2z^2\\
      \rule{1cm}{0cm} +15748y^4z^2-4374x^2z^4-4050y^2z^4-6568z^6),\\
      h:=xy(x+y)(x-y)(1819x^4+267x^2y^2+1819y^4+2404x^2z^2+2404y^2z^2+z^4).
   \end{array}
   \end{equation*}
   Since $I$ has three generators of the same degree, the containment criterion
   of Seceleanu \cite[Theorem 3.1]{Sec15} applies. As the details are rather
   tedious we postpone them to the forthcoming note \cite{MalSzp}.

\paragraph*{Acknowledgement.}
   We would like to thank Prof. Michael Cuntz for introducing
   us to the circle of ideas revolving around simplicial arrangements
   and for sharing with us a database of such arrangement. This
   research has began while we were visiting the University of Hannover in Spring 2017.
   It is the pleasure to thank M. Cuntz and P. Pokora for hospitality.
   We would like also to thank T. Szemberg for helpful remarks.\\
   Research of Malara was partially supported by National Science Centre, Poland, grant
   2016/21/N/ST1/01491.
   Research of Szpond was partially supported by National Science Centre, Poland, grant
   2014/15/B/ST1/02197.


\bigskip \small

\bigskip
   Grzegorz Malara,
   Department of Mathematics, Pedagogical University of Cracow,
   Podchor\c a\.zych 2,
   PL-30-084 Krak\'ow, Poland

\nopagebreak
   \textit{E-mail address:} \texttt{grzegorzmalara@gmail.com}

\bigskip
   Justyna Szpond,
   Department of Mathematics, Pedagogical University of Cracow,
   Podchor\c a\.zych 2,
   PL-30-084 Krak\'ow, Poland

\nopagebreak
   \textit{E-mail address:} \texttt{szpond@up.krakow.pl}


\end{document}